\documentclass[a4paper]{amsart}
\usepackage{graphicx}
\usepackage{amssymb}
\usepackage{amsmath}
\usepackage{amsthm}
\usepackage{amscd}
\usepackage[all,2cell]{xy}

\UseAllTwocells \SilentMatrices
\newtheorem{thm}{Theorem}[section]

\newtheorem{cor}[thm]{Corollary}
\newtheorem{lem}[thm]{Lemma}

\newtheorem{prop}[thm]{Proposition}
\newtheorem{prop-def}[thm]{Proposition-Definition}
\theoremstyle{definition}
\newtheorem{Def}[thm]{Definition}
\theoremstyle{remark}
\newtheorem{rmk}[thm]{\bf Remark}
\newtheorem{ex}[thm]{\bf Example}
\numberwithin{equation}{section}

\def\ot{\otimes}

\def\1{\textbf{1}}

\begin{document}
\title[$k$-torsionfree modules and Frobenius extensions]
{$k$-torsionfree modules and Frobenius extensions}
\author [Zhibing Zhao]
{Zhibing Zhao}


\thanks{Center for Pure Mathematics, School of Mathematical Sciences, Anhui University, Hefei 230601, Anhui, PR China}
 \thanks{E-mail:zbzhao$\symbol{64}$ahu.edu.cn (Z. B. Zhao).}
\maketitle

\dedicatory{}%
\commby{}%
\begin{abstract} Let $R/S$ be a Frobenius extension and $k$ be a positive integer. We prove that an $R$-module is $k$-torsionfree  if and only if so is its underlying $S$-module. As an application, we obtain that if $S$ is a quasi $k$-Gorenstein ring then so is $R$, but the converse does not hold in general.
\end{abstract}

\section{Introduction}

The notion of $k$-torsionfree modules was firstly introduced by Auslander and Bridger in \cite{AB}, which plays a crucial role in their theory on the stable category of finitely generated modules over a two-sided Noetherian ring.  It is well-known that, as the origin of Gorenstein homological algebra, the notion of modules of G-dimension zero was defined in terms of $k$-torsionfreeness \cite{AB}.

Let $R$ be a two-sided Noetherian ring and $R$-${\rm mod}$ the category of finitely generated left $R$-modules. For a module $M\in R$-${\rm mod}$, there is a projective resolution $\xymatrix@C=0.5cm{P_1\ar[r]^{f} & P_0\ar[r] & M  \ar[r] & 0}$ in $R$-${\rm mod}$. Putting $(-)^*={\rm Hom}_R(-,R)$, we have an exact sequence $\xymatrix@C=0.5cm{0 \ar[r] & {M^*}\ar[r] &P_0^*\ar[r]^{f^*}&P_1^* \ar[r] &{\rm Coker} {f^* }\ar[r] & 0}$, and call ${\rm Coker} {f^*}$ the transpose of $M$, denote it by ${\rm Tr}M$. The $R$-module $M$ is said to be \emph{$k$-torsionfree} if ${\rm Ext}^i_{R^{\rm op}}({\rm Tr}M,R)=0$ for $1\leq i\leq k$.
 We denote the full subcategory of $R$-${\rm mod}$ consisting of all $k$-torsionfree modules by $\mathcal{T}^k(R)$. It is easy to see that $\mathcal{T}^k(R)$ is closed under finite direct sums and  direct summands.

We are concerned with the extension-closedness of the subcategory $\mathcal{T}^k(R)$. Recall that a subcategory $\mathcal{X}$ of $R$-${\rm mod}$ is called \emph{extension-closed}, if the middle term $X$ of any short exact sequence $0\rightarrow X'\rightarrow X\rightarrow X''\rightarrow 0$ is in $\mathcal{X}$, provided the end terms $X', X''$ are in $\mathcal{X}$. Due to \cite[Theorem 4.7]{AR1} and \cite[Theorem 3.3]{Hu1}, the extension-closedness of the subcategory $\mathcal{T}^i(R)$ for $1\leq i\leq k$ equivalent to $R$ is right quasi $k$-Gorenstein.
  As a generalization of $k$-Gorenstein rings in the sense of Auslander, the notion of quasi $k$-Gorenstein rings was introduced by Huang in \cite{Hu2}.  We call a ring $R$ is  \emph{right  quasi $k$-Gorenstein} if the right  flat dimension of the $i$-st term in a minimal injective resolution of $R$ as a right $R$-module is less than or equal to $i$ for any $1\leq i\leq k$. The $k$-Gorenstein rings and their generalizations play important role in homological algebra, representation theory of algebras and non-commutative algebraic geometry \cite{AR, AR1, Hu2, HI}.

As a generalization of Frobenius algebras, the notion of  Frobenius extensions was firstly introduced by Kasch in \cite{Kas1}. The fundamental example of Frobenius extensions is the group
 algebras induced by a finite index subgroup. Frobenius algebras and extensions are of broad interest in many different areas, for example, they provide some connection between representation theory and knot theory \cite{Kad0}, and are used in study of Calabi-Yau properties of Cherednik algebras and quantum algebras \cite{BGS}.
  We refer to  \cite{Kad} for more details.

  It is well-known that many homological properties, such as the Gorenstein projectivity of modules, Iwanaga-Gorensteinness and the property of  FTF of rings and so on, are invariant under Frobenius extensions \cite{CR, GT, R, Z}. Gorenstein projective modules are central in Gorenstein homological algebra. A module is Gorenstein projective if and only if it is a module of G-dimension zero when it is a finitely generated module over a two-side Noetherian ring \cite{Chris}. And it is known that the notion of modules of G-dimension zero was defined in terms of $k$-torsionfreeness.  So there is a natural question: are the $k$-torsionfreeness of modules preserved under Frobenius extensions? We will investigate invariance of $k$-torsionfreeness of modules under Frobenius extensions and  give an affirmative answer of this question in this paper.

\noindent{\bf Setup and notation.} Throughout the paper, $R$ and $S$ are two-sided Noetherian rings and all modules are left modules if not specified otherwise. We denote by $R$-${\rm mod}$ (resp. ${\rm mod}$-$R$ ) the category of finitely generated left (resp. right) $R$-modules. We use pd$_R(M)$ (resp. id$_R(M)$) to denote the projective dimension (resp.  injective dimension) of $_RM$. For two right $R$-modules $A$ and $B$, we denote by ${\rm Hom}_{R^{\rm op}}(A,B)$ abelian group consisting of all right $R$-homomorphisms between them, where $R^{\rm op}$ is the opposite ring of $R$.

The following is the main result in this paper; see Theorem \ref{equivalence II}.

\noindent{\bf Theorem A.} Let $R/S$  be a Frobenius extension and $M$ an $R$-module and $k$ a positive integer. Then $M$ is $k$-torsionfree as an $R$-module if and only if $M$ is $k$-torsionfree as an $S$-module.

 Using Theorem A, we can prove that if $\mathcal{T}^i(S)$ is extension-closed then so is $\mathcal{T}^i(R)$ when $R/S$  is  a Frobenius extension, where $i$ is a positive integer. Note that a ring $R$ is right quasi $k$-Gorenstein if and only if $\mathcal{T}^i(R)$ is extension-closed for $1\leq i\leq k$. Thus, we obtain the following result, which maybe gives a new way to generate some new quasi $k$-Gorenstein rings (algebras); see Theorem \ref{quasi k-Gorenstein}.

\noindent{\bf Theorem B.} Let $R/S$  be a Frobenius extension and $k$ a positive integer. If $S$ is quasi $k$-Gorenstein, then so is $R$.

This paper is organized as follows. In Section 2, we give some notations  and preliminaries. We will investigate the transfer of $k$-torsionfreeness of modules  in Section 3, Theorem A are obtained; see Theorem \ref{equivalence II}. In Section 4,
we give an application of the main theorem A. We get  that if the base ring is quasi $k$-Gorenstein then so is the extension ring under Frobenius extensions; see Theorem \ref{quasi k-Gorenstein}. It maybe gives a new way to generate some new quasi $k$-Gorenstein rings (algebras). Moreover, we give a counterexample which shows that the converse of Theorem \ref{quasi k-Gorenstein} does not hold in general.  We discuss the transfer of $k$-Gorenstein property of rings under Frobenius extensions in Section 5. In Section 6, as another application of Theorem A, we give a new proof of the following known result: a module is of G-dimension zero over the extension ring if and only if so is its underlying module over the base ring under Frobenius extensions; see Theorem \ref{equivalent of G-dimension}.

\section{Preliminaries}

In this section, we recall some notations and collect some fundamental results which are used in this paper.

A {\it ring extension} $R/S$ is a ring homomorphism $\xymatrix@C=0.5cm{
   S \ar[r]^{l} &  R }$, which is required to send the unit to the unit. Thus $S$ is the base ring and $R$ is the extension ring. A ring extension is an algebra if $S$ is commutative and $l$ factor $S\rightarrow C(R)\hookrightarrow R$ where
   $C(R)$ is the center of $R$.
    The natural bimodule structure of ${_SR_S}$ is given by $s\cdot r\cdot s^\prime:=l(s)\cdot r\cdot l(s^\prime)$.
   Similarly, we can get some other module structures, for example $R_S$, ${_SR_R}$ and ${_RR_S}$, etc.

 Let $R/S$ be a ring extension and $M$ an $R$-module. Then $M$ is naturally an $S$-module, which is referred as the \emph{underlying} $S$-module. There is a natural surjective map $\pi: R\ot_SM\rightarrow M$ via $r\ot m\mapsto rm $ for any $r\in R$ and $m\in M$. It is easy check that $\pi$ is split as an $S$-homomorphism.  However, $\pi$ is not split as an $R$-homomorphism in general.

For a ring extension $R/S$, there is a {\it restriction functor} $Res: R$-${\rm Mod}\rightarrow S$-${\rm Mod}$ sending ${_RM}\mapsto {_SM}$, given by $s\cdot m:=l(s)\cdot m$. In the opposite direction, there are two natural functors as follows:

$(1)$ $\mathbb{T}={_RR}\ot_S- :S$-${\rm Mod}\rightarrow R$-${\rm Mod}$ is given by ${_SM}\mapsto {_RR}\ot_SM$.

$(2)$ $\mathbb{H}={\rm Hom}_S({_SR_R},-):S$-${\rm Mod}\rightarrow R$-${\rm Mod}$ is given by ${_SM}\mapsto {\rm Hom}_S({_SR_R},{_SM})$.

It is easy to check that both $(\mathbb{T},Res)$ and $(Res,\mathbb{H})$ are adjoint pairs.

Recall from \cite[Definition 2.1]{Kad} that a ring extension $R/S$ is said to be a \emph{Frobenius extension} if the functors $\mathbb{T}$ and $\mathbb{H}$ above are naturally equvilent. It is well-known that the condition is equivalent to the conditions that ${_SR}$ is finitely generated projective and ${_RR_S}\cong {\rm Hom}_S({_SR_R},{_SS_S})$. By \cite[Theorem 1.2]{Kad}, these conditions are also equivalent to the corresponding properties on the opposite sides.

There are other examples of Frobenius extensions including Hopf subalgebras, finite extensions of
enveloping algebras of Lie super-algebra and finite extensions of enveloping algebras of Lie coloralgebras  etc \cite{FMS,Sch}. More examples of Frobenius extensions can be found in \cite[Example 2.4]{Z}.

As a generalization of the classical torsionless  and reflexive modules, the notion of $k$-torsionfree modules was firstly introduced by Auslander and Bridger in \cite{AB}. Due to Auslander, there is an exact sequence
$$\xymatrix@C=0.5cm{0 \ar[r] & {\rm Ext}^1_{R^{\rm op}}({\rm Tr}M, R)\ar[r] & M\ar[r]^{\sigma_M}& M^{**}\ar[r] &{\rm Ext}^2_{R^{\rm op}}({\rm Tr}M, R)\ar[r] & 0}$$
 for any module $M$ in $R$-${\rm mod}$, where $\sigma_M:M\rightarrow M^{**}$ is the evaluation map via $\sigma_M(m)(f)=f(m)$ for any $f\in M^*$ and $m\in M$ (see \cite[Propostion 6.3]{A}). Then $M$ is 1-torsionfree  if and only if it is torsionless, and $M$ is 2-torsionfree if and only if it is reflexive. In \cite{AB}, Auslander and Bridger have obtained a description of $k$-torsionfreeness in terms of approximations.

The following definition was introduced firstly by Auslander and Smal${\o}$ in \cite{AS}.
\begin{Def}{\label{approximation}} Let $\mathcal{C}$ be a full subcategory of $R$-${\rm mod}$ and $C\in \mathcal{C}$, $M\in R$-${\rm mod}$. An $R$-homomorphism $M\rightarrow C$ is said to be a \emph{left $\mathcal{C}$-approximation} of $M$ if ${\rm Hom}_R(C, X)\rightarrow{\rm Hom}_R(M, X)$ is epic for any $X\in\mathcal{C}$. An subcategory $\mathcal{C}$ is said to be \emph{covariantly finite} in $R$-${\rm mod}$ if every module in $R$-${\rm mod}$ has a left $\mathcal{C}$-approximation. Dually, An $R$-homomorphism $C\rightarrow M$ is said to be a \emph{right $\mathcal{C}$-approximation} of $M$ if ${\rm Hom}_R(X, C)\rightarrow{\rm Hom}_R(X, M)$ is epic for any $X\in\mathcal{C}$. An subcategory $\mathcal{C}$ is said to be \emph{contravariantly finite} in $R$-${\rm mod}$ if every module in $R$-${\rm mod}$ has a right $\mathcal{C}$-approximation. An subcategory $\mathcal{C}$ is said to be \emph{functorially  finite} in $R$-${\rm mod}$ if it is both covariantly finite and contraviantly finite in $R$-${\rm mod}$.
\end{Def}

\begin{lem}{\rm(}{\cite[Theorem 2.17]{AB}}\rm{)}{\label{equivalence of k-tosionfree}} For a module $M\in R$-${\rm mod}$ and a positive integer $k$, the following statements are equivalent.

$(1)$ $M$ is $k$-torsionfree.

$(2)$ There exists an exact sequence $\xymatrix@C=0.5cm{0 \ar[r] & M\ar[r]^{f_1} &{P_1}\ar[r]^{f_2~} &  \cdots \ar[r]^{f_k} &{P_k}}$ with $P_i$ finitely generated projective, such that ${\rm Im} f_i \rightarrow P_i$ is a left ${\rm add}R$-approximation of ${\rm Im} f_i$ for $1\leq i\leq k$, where ${\rm add}R$ is the subcategory of $R$-${\rm mod}$ that consisting of all modules isomorphic to direct summands of finite direct sum of copies of ${_RR}$.

$(3)$ There exists an exact sequence $\xymatrix@C=0.5cm{0 \ar[r] & M\ar[r]^{f_1} &{F_1}\ar[r]^{f_2~} &  \cdots \ar[r]^{f_k} &{F_k}}$ with $F_i$ finitely generated free, such that ${\rm Im} f_i \rightarrow F_i$ is a left ${\rm add}R$-approximation of ${\rm Im} f_i$ for $1\leq i\leq k$.
\end{lem}

This description of $k$-torsionfree modules lead to the related notion of $k$-syzygy modules. A module $M$ is called $k$-syzygy if there exists an exact sequence $0\rightarrow M\rightarrow P_{k-1}\rightarrow\cdots\rightarrow P_0\rightarrow N \rightarrow 0$ in $R$-${\rm mod}$ with $P_i$ projective for $0\leq i\leq k-1$. We denote by $\Omega^k(R)$ the full subcategory of $R$-${\rm mod}$ whose objects are the $k$-syzygy modules in $R$-${\rm mod}$. It is clear that  $\mathcal{T}^k(R)\subseteq\Omega^k(R)$ by the lemma above. But a $k$-syzygy module is not $k$-torsionfree in general.

It is well-known that $\Omega^k(R)$ is a covariantly finite subcategory   of $R$-${\rm mod}$; see \cite[Corollary 1.8]{AR2}. Covariantly finite subcategories of $R$-${\rm mod}$ have been proved to be of interest in various situations, for instance, in the representation of Artin algebras and commutative ring theory. They are particularly interesting when they are extension-closed. If $R$ is an Artin algebra and $\Omega^i(R)$ is closed under extensions, then $\Omega^i(R)$ is functorially finite resolving, and
$\Omega^i(R)$ has almost split sequences; see \cite{AR}. It follows form \cite[Theorem 3.3]{Hu1} that $\mathcal{T}^i(R)$ is extension-closed for $1\leq i\leq k$ if and only if so is $\Omega^i(R)$ for $1\leq i\leq k$, in this case,
$\mathcal{T}^i(R)=\Omega^i(R)$ and $R$ is called a right quasi $k$-Gorenstein ring. This is the main reason why we are concerned with extension-closedness of $\mathcal{T}^k(R)$.

The following observation gives a recurrence relation of $k$-torsionfreeness.
\begin{prop}\label{recurrence of k-torsionfree} Let $\xymatrix@C=0.5cm{0\ar[r] & M\ar[r]^{f} &P\ar[r] & N\ar[r] &0}$ be a short exact sequence in $R$-${\rm mod}$ with $f^*$  epic and let $k\geq 1$ be an integer.  Then $M$ is $k$-torsionfree if and only if $N$ is $(k-1)$-torsionfree.
\end{prop}
\noindent{\bf Proof.} $(\Rightarrow)$ If $k=1$, it is trivial.

For the case of  $k=2$, we will show that if $M$ is reflexive then $N$ is torsionless. Since $f^*$ is epimorphic, applying by $(-)^*={\rm Hom}_R(-,R)$, we get an exact sequence
$$\xymatrix@C=0.5cm{0\ar[r] & N^*\ar[r] &P^*\ar[r]^{f^*} & M^*\ar[r] &0}.$$

Consider the following commutative diagram with exact rows
$$\xymatrix{
 0\ar[r]& M\ar[d]_{\sigma_M} \ar[r]^{f} & P\ar[d]_{\sigma_P} \ar[r] &N\ar[d]_{\sigma_N}\ar[r]& 0\\
  0\ar[r]& M^{**} \ar[r] &  P^{**}\ar[r] & N^{**}},$$
by Snake Lemma,  $\sigma_N$ is a monomorphism since $\sigma_M$ and $\sigma_P$ are isomorphisms. Hence $N$ is a torsionless module.

Now, we suppose that $k\geq 3$. Since $M$ is a $k$-torsionfree module, $M$ must be 2-torsionfree. By the proof of case for $k=2$, $N$ is 1-torsionfree, and so ${\rm Ext}^{1}_{R^{op}}({\rm Tr}N, R)=0$.  Let $\xymatrix@C=0.5cm{P_1\ar[r] & P_0\ar[r] & M \ar[r] & 0}$ be a projective resolution of $M$, then we have an exact sequence
$\xymatrix@C=0.5cm{0 \ar[r] & {M^*}\ar[r] &P_0^*\ar[r]&P_1^* \ar[r] & {\rm Tr}M\ar[r] & 0}$ in mod-$R$. By dimension shifting and assumption, we get
${\rm Ext}^i_{R^{\rm op}}(M^*,R)\cong{\rm Ext}^{i+2}_{R^{\rm op}}({\rm Tr}M,$\\$R)=0$
for $1\leq i\leq k-2$. Similarly, we get ${\rm Ext}^i_{R^{\rm op}}(N^*, R)\cong{\rm Ext}^{i+2}_{R^{\rm op}}({\rm Tr}N, R)$ for $i\geq 1$.

Since $f^*$ is an epimorphism, the sequence $\xymatrix@C=0.5cm{0\ar[r] & N^*\ar[r] &P^*\ar[r]^{f^*} & M^*\ar[r] &0}$  is exact in mod-$R$. Then we have ${\rm Ext}^i_{R^{\rm op}}(N^*, R)\cong{\rm Ext}^{i+1}_{R^{\rm op}}(M^*, R)$ for $i\geq 1$. Hence,
${\rm Ext}^{i}_{R^{\rm op}}({\rm Tr}N, R)\cong{\rm Ext}^{i-2}_{R^{\rm op}}(N^*, R)\cong{\rm Ext}^{i-1}_{R^{\rm op}}(M^*, R)\cong{\rm Ext}^{i+1}_{R^{\rm op}}({\rm Tr}M, R)=0$ for $2\leq i\leq k-1$. So we have ${\rm Ext}^{i}_{R^{\rm op}}({\rm Tr}N, R)=0$ for $1\leq i\leq k-1$. That is, $N$ is a $(k-1)$-torsionfree module.

$(\Leftarrow)$ If $N$ is $(k-1)$-torsionfree, then there exists an exact sequence $\xymatrix@C=0.5cm{0 \ar[r] & N\ar[r]^{f_1} &{P_1}}$\\$\xymatrix@C=0.5cm{\ar[r]^{f_2~} &  \cdots \ar[r]^{f_{k-1}} &{P_{k-1}}}$ with each $P_i$ finitely generated projective, such that ${\rm Im} f_i \rightarrow P_i$ is a left ${\rm add}R$-approximation of ${\rm Im} f_i$ for $1\leq i\leq k-1$ by Lemma \ref{equivalence of k-tosionfree}. Consider the exact sequence $\xymatrix@C=0.5cm{0 \ar[r] & M \ar[r]^f & P\ar[r]^{f_1} &{P_1}\ar[r]^{f_2~} &  \cdots \ar[r]^{f_{k-1}} &{P_{k-1}}}$, $M(\cong{\rm Im}f)\rightarrow P$ is a left ${\rm add}R$-approximation since $f^*$ is an epimorphism and ${\rm Im} f_i \rightarrow P_i$ is a left ${\rm add}R$-approximation of ${\rm Im} f_i$ for $1\leq i\leq k-1$. Hence $M$ is $k$-torsionfree by Lemma \ref{equivalence of k-tosionfree} again.  \hfill$\square$

\section{$k$-torsionfree modules}
In this section, we will show that $k$-torsionfreeness of a module is preserved under Frobenius extensions, that is, a module is $k$-torsionfree over the extension ring if and only if so is its underlying module over the base ring under Frobenius extensions.

The following observation is needed.
\begin{lem}{\label{lemma of torsionfree}} Let $M$ be in $R$-${\rm mod}$. Then $M$ is $k$-torsionfree if and only if there exists an exact sequence
$$\xymatrix@C=0.5cm{0 \ar[r] & M\ar[r]^{f_1} &{P_1}\ar[r]^{f_2~} &  \cdots \ar[r]^{f_k} &{P_k}\ar[r] & T_k \ar[r] & 0}$$
with $P_i$ finitely generated projective and  ${\rm Ext}^{j}_R(T_k, R)=0$ for $1\leq j\leq k$.
\end{lem}
\noindent{\bf Proof.} If $M$ is $k$-torsionfree, then there exists an exact sequence  $$\xymatrix@C=0.5cm{0 \ar[r] & M\ar[r]^{f_1} &{P_1}\ar[r]^{f_2~} &  \cdots \ar[r]^{f_k} &{P_k}\ar[r] & T_k \ar[r] & 0}$$
 with $P_i$ finitely generated projective, such that ${\rm Im} f_i \rightarrow P_i$ is a left ${\rm add}R$-approximation of ${\rm Im} f_i$ for $1\leq i\leq k$ by Lemma \ref{equivalence of k-tosionfree}. Taking ${\rm Coker}f_i=T_i$, we have that ${\rm Ext}^1_R(T_i, R)=0$ for $1\leq i\leq k$. By dimension shifting,  ${\rm Ext}^j_R(T_i, R)\cong{\rm Ext}^{k-i+j}_R(T_k, R)$, and we obtain that ${\rm Ext}^{j}_R(T_k, R)=0$ for $1\leq j\leq k$.

Conversely, if there exists an exact sequence
$\xymatrix@C=0.5cm{0 \ar[r] & M\ar[r]^{f_1} &{P_1}\ar[r]^{f_2~} & \cdots \ar[r]^{f_k} &{P_k}\ar[r] & T_k}$
$\xymatrix@C=0.5cm{\ar[r] & 0}$
with $P_i$ finitely generated projective and  ${\rm Ext}^{j}_R(T_k, R)=0$ for $1\leq j\leq k$, then ${\rm Ext}^1_R(T_i, R)\cong{\rm Ext}^{k-i+1}_R(T_k, R)=0$
for $1\leq i\leq k$ where $T_i={\rm Coker}f_i$. It follows that ${\rm Im} f_i \rightarrow P_i$ is a left ${\rm add}R$-approximation of ${\rm Im} f_i$ for $1\leq i\leq k$, and $M$ is $k$-torsionfree by Lemma \ref{equivalence of k-tosionfree} again.\hfill$\square$

\begin{prop}{\label{necessity of main thm}} Let $R/S$  be a Frobenius extension and $M$ an $R$-module and $k$ a positive integer. If  $M$ is $k$-torsionfree as an $R$-module, then  $M$ is also $k$-torsionfree as the underlying $S$-module.
\end{prop}
\noindent{\bf Proof.} Suppose that $M$ is a $k$-torsionfree $R$-module. By Lemma \ref{lemma of torsionfree}, there exists an exact sequence
\begin{equation}{\label{3.1}}
\xymatrix@C=0.5cm{0 \ar[r] & {_RM}\ar[r]^{f_1} &{P_1}\ar[r]^{f_2~} &  \cdots \ar[r]^{f_k} &{P_k}\ar[r] & T_k \ar[r] & 0}
\end{equation}
in $R$-${\rm mod}$ with $P_i$ finitely generated projective and
${\rm Ext}^{j}_R(T_k, R)=0$ for $1\leq j\leq k$. Applying to (\ref{3.1}) by restriction functor $_SR\ot_R-$, we get the following exact sequence
\begin{equation}
\xymatrix@C=0.5cm{0 \ar[r] & {_SM}\ar[r]^{f_1} &{P_1}\ar[r]^{f_2~} &  \cdots \ar[r]^{f_k} &{P_k}\ar[r] & T_k \ar[r] & 0}
\end{equation}
in $S$-${\rm mod}$ with $P_i$ finitely generated projective and $T_k={\rm Coker}f_k$. Since $R/S$  is a Frobenius extension,  ${\rm Hom}_S({_SR_R},-)\cong{_RR\ot_S-}$. By the adjoint isomorphism, for $1\leq j\leq k$, we have
\begin{align*}{\rm Ext}_S^j(T_k,S)
&\cong {\rm Ext}_S^j({_SR\ot_RT_k},S)\\
&\cong {\rm Ext}_R^j(T_k,{\rm Hom}_S({_SR_R},S))\\
&\cong {\rm Ext}_R^j(T_k,{_RR\ot_SS})\\
&\cong {\rm Ext}_R^j(T_k,{_RR})\\
&=0.
\end{align*}
Thus, $M$ is a $k$-torsionfree $S$-module by Lemma \ref{lemma of torsionfree}.\hfill$\square$

\begin{prop}{\label{equivalence I}} Let $R/S$  be a Frobenius extension and $M$ an $S$-module and $k$ a positive integer. If $M$ is $k$-torsionfree as an $S$-module then ${_RR}\ot_SM$ is $k$-torsionfree as an $R$-module.
\end{prop}
\noindent{\bf Proof.} Suppose that ${_SM}$ is a $k$-torsionfree $S$-module, then there exists an exact sequence
$$\xymatrix@C=0.5cm{0 \ar[r] & {_SM}\ar[r]^{f_1} &{P_1}\ar[r]^{f_2~} &{P_2}\ar[r] &\cdots \ar[r]^{f_k} &{P_k}\ar[r]& {T_k}\ar[r]& 0}$$
in $S$-${\rm mod}$ with $P_i$ finitely generated projective and ${\rm Ext}^i_S(T_k, S)=0$ for $1\leq i\leq k$ by Lemma \ref{lemma of torsionfree}.
Since $R/S$  is a Frobenius extension, $R_S$ is finitely generated projective, and so ${_RR\ot_S-}$ is an exact functor. We get the following exact sequence
$$\xymatrix@C=0.5cm{0 \ar[r] & {{_RR\ot_S}M}\ar[r]^{{R\ot_S}f_1} &{{_RR\ot_S}P_1}\ar[r]^{{R\ot_S}f_2~} & {{_RR\ot_S}P_2}\ar[r] & \cdots \ar[r]^{{R\ot_S}f_k} &{{_RR\ot_S}P_k}\ar[r]& {{_RR\ot_S}T_k}\ar[r]& 0}$$
in $R$-${\rm mod}$ with ${_RR\ot_S}P_i$ finitely generated projective for $1\leq i\leq k$ and ${{_RR\ot_S}T_k}\cong {\rm Coker}({R\ot_S}f_k)$.
By the adjoint isomorphism, for any $i\geq 0$, we have
\begin{align*}{\rm Ext}^i_R({_RR\ot_S}T_k,R)
&\cong{\rm Ext}^{i}_S({_ST_k}, {\rm Hom}_R({_RR_S}, R))\\
&\cong {\rm Ext}^i_S({_ST_k}, {_SR}).
\end{align*}
Since ${_SR}$ is finitely generated projective as an $S$-module, $0={\rm Ext}^i_S({_ST_k}, {_SR})\cong{\rm Ext}^i_R({_RR\ot_S}T_k,R)$ for $1\leq i\leq k$.
Then ${_RR\ot_S}M$ is $k$-torsionfree as an $R$-module by Lemma \ref{lemma of torsionfree}.
\hfill$\square$

\begin{rmk} We mention that the converse of the proposition above does not hold in general. Let $R$ and $R^\prime$ be  two-sided Noetherian rings. An module $(M, M^\prime)$ over the direct product ring $R\times R^\prime$ is $k$-torsionfree if and only if so are both $M$ and $M^\prime$ over $R$ and $R^\prime$, respectively.

Consider the canonical projection ${\rm Pr}: (R, R^\prime)\rightarrow R$. That is a Frobenius extension from $R\times R^\prime$ to $R$. It is clear that $M$ is $k$-torsionfree as an $R$-module can not induced that $(M, M^\prime)$ is $k$-torsionfree as an $R\times R^\prime$-module.
\end{rmk}

The following is the main theorem in this section, which shows that the $k$-torsionfreeness of modules is preserved under Frobenius extensions.
\begin{thm}{\label{equivalence II}} Let $R/S$  be a Frobenius extension and $M$ an $R$-module and $k$ a positive integer. Then $M$ is $k$-torsionfree as an $R$-module if and only if $M$ is $k$-torsionfree as the underlying $S$-module.
\end{thm}
\noindent{\bf Proof.} $(\Rightarrow)$ By Proposition \ref{necessity of main thm}.

$(\Leftarrow)$ If $k=1$, $M$ is a torsionless $S$-module, which is equivalent to $M$ is cogenerated by $_SS$. Thus, there exists an exact sequence
$\xymatrix@C=0.5cm{0 \ar[r] & M\ar[r]^{f} &S^{I}}$ with some indexed set $I$. Applying by the functor ${\rm Hom}_S({_SR_R},-)$, since ${_SR}$ is finitely generated projective, we have an exact sequence
 $\xymatrix@C=0.5cm{0 \xrightarrow{}{\rm Hom}_S({R_R},M)\xrightarrow{{\rm Hom}_S({R},f)}{{\rm Hom}_S({R_R},S^I})}$ in $R$-${\rm mod}$. By assumption, ${_RR\ot_S-}\cong{\rm Hom}_S({_SR_R},-)$, we have
\begin{align*}
{\rm Hom}_S({_SR_R},S^I)
&\cong\prod\limits_{i\in I}{\rm Hom}_S({_SR_R},S)\\
&\cong\prod\limits_{i\in I}{_RR\ot_SS}\\
&\cong{_RR}^I.
\end{align*}

Because there is an $R$-monomorphism $i:{_RM}\rightarrow {\rm Hom}_S({_SR_R},{_SM})$ via $i(m)(r)=rm$ for any $m\in M$ and $r\in R$, we have an exact sequence
$0 \rightarrow{_RM}\xrightarrow{{\rm Hom}_S({R},f)\cdot i}{{\rm Hom}_S({_SR_R},S^I})\cong{_RR^{I}}$.  Hence ${_RM}$ is cogenerated by $_RR$, and so $M$ is a torsionless $R$-module.

Now we assume that $k>1$, $M$ is a $k$-torsionfree module as an $S$-module, there exists an exact sequence
$\xymatrix@C=0.5cm{0 \xrightarrow{}{_SM}\xrightarrow{f_1}{P_1}\xrightarrow{f_2~}{P_2}\xrightarrow{}\cdots}\xrightarrow{f_k}{P_k}$
in $S$-${\rm mod}$ with $P_i$ finitely generated projective and ${\rm Im}f_i\hookrightarrow P_i$ is a left ${\rm add}S$-approximation of ${\rm Im}f_i$ for $1\leq i\leq k$ by Lemma \ref{equivalence of k-tosionfree}.
Applying by the functor ${_RR\ot_S-}$, we get the following exact sequence
$$\xymatrix@C=0.8cm{0 \xrightarrow{}{{_RR\ot_S}M}\xrightarrow{{R\ot_S}f_1}{{_RR\ot_S}P_1}\xrightarrow{{R\ot_S}f_2~}{{_RR\ot_S}P_2}\xrightarrow{} \cdots\xrightarrow{{R\ot_S}f_k}{{_RR\ot_S}P_k}}$$
in $R$-${\rm mod}$ with ${_RR\ot_S}P_i$ finitely generated projective for $1\leq i\leq k$.
By Proposition \ref{equivalence I} and its proof, ${_RR\ot_S}M$ is $k$-torsionfree as an $R$-module and ${\rm Im}({R\ot_S}f_i)\hookrightarrow {{_RR\ot_S}P_i}$ is a left ${\rm add}R$-approximation of ${\rm Im}({R\ot_S}f_i)$ for $1\leq i\leq k$.

Note that there is an $R$-monomorphism $i:{_RM}\rightarrow {\rm Hom}_S({_SR_R},{_SM})$ via $i(m)(r)=rm$ for any $m\in M$ and $r\in R$, which is split when we restrict it as an $S$-homomorphism. Meanwhile, there is an $R$-epimorphism $\pi:{_RR\ot_SM}\rightarrow {_RM}$ via $\pi(r\ot m)=rm$ for any $m\in M$ and $r\in R$, which is also split when we restrict is as an $S$-homomorphism.
Since ${_RR\ot_SM} \cong {\rm Hom}_S({_SR_R},M)$, there is an $R$-monomorphism, which still denoted it by $i$, $i:{_RM}\hookrightarrow {_RR\ot_SM}$ and it is split
as a homomorphism of $S$-modules.
Denote ${_RR\ot_S}f_1$ by $f$, then there exists an exact sequence
\begin{equation}{\label{3.3}}
\xymatrix@C=0.5cm{0 \xrightarrow{} {_RM}\xrightarrow{f\cdot i}{_RR\ot_SP_1}\xrightarrow{}T\xrightarrow{}0},
\end{equation}
where $T\cong {\rm Coker}(f\cdot i)$.
 We claim that $\xymatrix@C=0.5cm{0 \xrightarrow{}{_RM}\xrightarrow{f\cdot i}{_RR\ot_SP_1}}$ is a left ${\rm add}R$-approximation of $_RM$.

  Applying to (\ref{3.3}) by restriction functor, we get an exact sequence
  $$\xymatrix@C=0.5cm{0 \xrightarrow{} {_SM}\xrightarrow{f\cdot i}{_SR\ot_SP_1}\xrightarrow{}T\xrightarrow{}0}$$ in $S$-${\rm mod}$. Let $h:M\rightarrow S$ be an $S$-homomorphism. Then $h\cdot\pi$ is an $S$-homomorphism from ${_SR\ot_SM}$ to $S$, where $\pi$ is the split $S$-homomorphism from ${R\ot_SM}$ to $M$. By necessity, ${_SR\ot_SM}$ is $k$-torsionfree and $f$ (as an $S$-homomorphism) is a left ${\rm add}S$-approximation of ${_SR\ot_SM}$. Hence,
 there exists an $S$-homomorphism $g:{_SR\ot_SP_1}\rightarrow S$ such that $g\cdot f=h\cdot\pi$
 $$\xymatrix{
  & S  &  \\
  0\ar[r]& {_SM} \ar[u]_{\forall h}\ar[r]^{f\cdot i} &  {_SR\ot_SP_1}\ar@{..>}[ul]_{\exists g} \ar[r] & {_ST} \ar[r] & 0 \\
  &{_SR\ot_SM} \ar[u]_{\pi} \ar[ur]_{f} &  }.$$
Then, $g\cdot(f\cdot i)=(g\cdot f)\cdot i=(h\cdot\pi)\cdot i=h\cdot(\pi\cdot i)=h$. This means that $f\cdot i:{_SM}\hookrightarrow{_SR\ot_SP_1}$ is a left ${\rm add}S$-approximation of $_SM$.
Since $_SM$ is $k$-torsionfree, $T$ is $(k-1)$-torsionfree as an $S$-module by Proposition \ref{recurrence of k-torsionfree}. By induction hypothesis, $T$ is also a $(k-1)$-torsionfree $R$-module. Since $f\cdot i:{_SM}\hookrightarrow{_SR\ot_SP_1}$ is a left ${\rm add}S$-approximation of $M$, we have
\begin{align*}0
&={\rm Ext}^1_S({_ST},S)\\
&\cong{\rm Ext}^{1}_S({_SR\ot_RT},S)\\
&\cong{\rm Ext}^{1}_R({_RT}, {\rm Hom}_S({_SR_R}, S))\\
&\cong {\rm Ext}^1_R({_RT}, {_RR\ot_S}S)\\
&\cong {\rm Ext}^1_R({_RT}, {_RR}).
\end{align*}
Applying by Hom$_R(-, R)$ to the exact sequence (3.3), we have ${\rm Hom}_R({_RR\ot_SP_1},R)\rightarrow{\rm Hom}_R({_RM},R)\rightarrow {\rm Ext}^1_R({_RT}, {_RR})=0$. Hence $f\cdot i:{_RM}\hookrightarrow{_RR\ot_SP_1}$ is a left ${\rm add}R$-approximation of ${_RM}$. Since $T$ is  a $(k-1)$-torsionfree $R$-module and the sequence (3.3) is exact,  ${_RM}$ is a $k$-torsionfree module as an $R$-module by Proposition \ref{recurrence of k-torsionfree}. We finish the proof. \hfill$\square$

For an $R$-module $M$, it is called an \emph{$\infty$-torsionfree} module if it is $k$-torsionfree for any positive integer $k$. By the theorem above,  we have the following.

\begin{cor}{\label{torsionfree equ II}} Let $R/S$ be a Frobenius extension and $M$ an $R$-module. Then we have that $M$ is $\infty$-torsionfree as an $R$-module if and only if $M$ is $\infty$-torsionfree as the underlying $S$-module.
\end{cor}

\section{quasi $k$-Gorenstein rings}
In this section, we will give an application of the results in Section 3. Recall that a ring $R$ is called \emph{right (left) quasi $k$-Gorenstein} if the right (left) flat dimension of the $i$-st term in a minimal injective resolution of $R$ as a right (left) $R$-module is less than or equal to $i$ for any $1\leq i\leq k$. If a ring $R$ is both left and right quasi $k$-Gorenstein, then $R$ is said to be a quasi $k$-Gorenstein ring. It is well-known that quasi $k$-Gorenstein rings are not left-right symmetric in general; see \cite [Example 2]{Hu2}.

The following is the main result in this section, which is possible to generate new left (right) quasi  $k$-Gorenstein rings.

\begin{thm}{\label{quasi k-Gorenstein}} Let $R/S$  be a Frobenius extension and $k$ a positive integer. If $S$ is a right (left) quasi $k$-Gorenstein ring, then so is $R$.
\end{thm}

The following result gives equivalent characterization of right quasi $k$-Gorenstein rings.

\begin{lem}${\label{extension closed of k-tosionfree}}($\cite[Theorem 4.7]{AR1} and \cite[Theorem 3.3]{Hu1}$)$ Let $R$ be a ring. The following statements are equivalent for a positive integer $k$.

$(1)$ $\Omega^i(R)$ is extension-closed for $1\leq i\leq k$;

$(2)$ ${\rm add}\Omega^i(R)$ is extension-closed for $1\leq i\leq k$;

$(3)$ If $0\rightarrow R\rightarrow I'_0\rightarrow I'_1\rightarrow\cdots\rightarrow I'_i\rightarrow\cdots$ is a minimal injective resolution of $R$ as a right $R$-module, then $r.{\rm fd}_R(I'_i)\leq i+1$ for $0\leq i\leq k-1$;

$(4)$ $\mathcal{T}^i(R)$ is extension-closed for $1\leq i\leq k$.

Furthermore, if one of equivalent conditions is satisfied, then $\mathcal{T}^i(R)=\Omega^i(R)={\rm add}\Omega^i(R)$ for $1\leq i\leq k$.
\end{lem}

By the lemma above, a ring $R$ is right quasi $k$-Gorenstein if and only if $\mathcal{T}^i(R)=\Omega^i(R)$ are extension-closed for $1\leq i\leq k$. This is the reason we concern the extension-closedness of the subcategory $\mathcal{T}^i(R)$. In fact, by symmetry of Frobenius extensions, Theorem \ref{quasi k-Gorenstein} is a consequence of the following proposition.

\begin{prop}{\label{equivalent of extension closure}}Let $R/S$  be a Frobenius extension and $k$ a positive integer. If $\mathcal{T}^i(S)$ is extension-closed for $1\leq i\leq k$, then so is $\mathcal{T}^i(R)$ for $1\leq i\leq k$.
\end{prop}
\noindent{\bf Proof.} Let $0\rightarrow A\rightarrow B\rightarrow C\rightarrow 0$ be an exact sequence in $R$-mod such that $A, C$ are $i$-torsionfree modules for $1\leq i\leq k$. We want show that $B$ is also $i$-torsionfree for $1\leq i\leq k$. Applying by the restriction  functor $_SR\ot_R-$, we have the  exact sequence
$0\rightarrow A\rightarrow B\rightarrow C\rightarrow 0$
in $S$-mod. It follows from Proposition \ref{necessity of main thm} that $A$ and $C$ are $i$-torsionfree as $S$-modules for $1\leq i\leq k$.
Since $\mathcal{T}^i(S)$ is extension-closed for $1\leq i\leq k$, $B$ is also $i$-torsionfree as an $S$-module for $1\leq i\leq k$. So $B$ is $i$-torsionfree as an $R$-module by Theorem \ref{equivalence II} for $1\leq i\leq k$.   \hfill$\square$

\begin{rmk} The converse of the proposition above does not hold in general. Let $R$ and $R^\prime$ be  two-sided Noetherian rings. Take $S$ to be the direct product ring $R\times R^\prime$ of $R$ and $R^\prime$.
Consider the canonical projection ${\rm Pr}: S=R\times R^\prime\rightarrow R$. It is a Frobenius extension from $S$ to $R$. For any $1\leq i\leq k$, it is clear that $\mathcal{T}^i(S)$ is extension-closed if and only if $\mathcal{T}^i(R)$ and $\mathcal{T}^i(R^\prime)$ are extension-closed, respectively. Therefore, the quasi $k$-Gorenstein property of rings does not reflect under Frobenius extensions in general.
\end{rmk}

By Lemma \ref{extension closed of k-tosionfree}, the following is a direct consequence  of the propostion above.
\begin{cor} Let $R/S$  be a Frobenius extension and $k$ a positive integer. If  $\Omega^i(S)$ is extension-closed for $1\leq i\leq k$ then so is $\Omega^i(R)$ for $1\leq i\leq k$.  Thus, we obtain that if  $\Omega^i(S)$ are functorially finite in $S$-${\rm mod}$ for $1\leq i\leq k$ then $\Omega^i(R)$ are functorially finite in $R$-${\rm mod}$ for $1\leq i\leq k$.
\end{cor}

Similar to the proof of \cite[Theorem 3.10]{FGR}, we may prove that a ring $R$ is quasi $k$-Gorenstein if and only if $\left(
\begin{array}{cc}
R & 0\\
R & R
\end{array}
\right)$ is quasi $k$-Gorenstein. Theorem \ref{quasi k-Gorenstein} gives a new method to generate new right (left) quasi $k$-Gorenstein rings.
We end this section with the following examples.

\begin{ex} $(1)$ Let $R$ be a two-sided Noetherian ring and $M_n(R)$ the $n\times n$ matrix ring over $R$.  Then $M_n(R)$ is a Frobenius extension of $R$. By Theorem \ref{quasi k-Gorenstein}, if  $R$ is right (left) quasi $k$-Gorenstein then so is $M_n(R)$. And $\Omega^i(M_n(R))$ is extension-closed for $1\leq i\leq k$ can be deduced from that the extension-closedness of  $\Omega^i(R)$ for $1\leq i\leq k$.

$(2)$ Let $G$ be a group and $H$ be subgroup with finte index. Let $R$ be a ring. Denote by $RG$ and $RH$ the group rings. Then $RH\hookrightarrow RG$ is a Frobenius extension.
 It follows from Theorem \ref{quasi k-Gorenstein} that  if $RH$ is right (left) quasi $k$-Gorenstein then so is $RH$.

$(3)$ Let $R$ be a two-sided Noetherian ring and $G$ a finite group. We have that the skew group ring $R*G$ is a Frobenius extension of $R$.
If $R$ is right (left) quasi $k$-Gorenstein then so is $R*G$.

$(4)$ Let $A$ be an Artin ring and $A[x]/(x^t)$ ($t\geq 2$) the quotient of the polynomial
ring, where $x$ is a variable which is supposed to commute with all the elements of $A$. Then the truncated polynomial extension $R\rightarrow A[x]/(x^t)$ is a Frobenius extension. We have that if $A$ is  right (left) quasi $k$-Gorenstein then so is $Q$.

$(5)$ Let $A$ be a central separable Artin algebra over center $C$. Then $A$ is a Frobenius extension of $C$, see \cite{Kad0}. By Theorem \ref{quasi k-Gorenstein}, if $C$ is  right (left) quasi $k$-Gorenstein then so is $A$.
\end{ex}

\section{$k$-Gorenstein rings}

In this section, we will investigate the transfer of $k$-Gorenstein property of rings in sense of Auslander under Frobenius extensions.

 Motivated by the theory of commutative noetherian Gorenstein rings, Auslander introduced the notion of $k$-Gorenstein rings \cite{FGR}. Let $0\rightarrow R\rightarrow I_0\rightarrow I_1\rightarrow\cdots\rightarrow I_i\rightarrow\cdots$ be a minimal injective resolution of $R$ as an $R$-module and
$0\rightarrow R\rightarrow {I_0^\prime}\rightarrow I_1^\prime\rightarrow\cdots\rightarrow I_i^\prime\rightarrow\cdots$ be a minimal injective resolution of $R$ as a right $R$-module. The ring $R$ is called a \emph{$k$-Gorenstein ring} if $r.{\rm fd}_{R^{{\rm op}}}(I_i^\prime)\leq i$ for any $0\leq i\leq k-1$, and $R$ is said to be an \emph{Auslander ring} if it is $k$-Gorenstein for all $k$. An important fact is that the notion of $k$-Gorenstein rings is left-right symmetric by \cite[Auslander's Theorem 3.7]{FGR}.

It is clear that a $k$-Gorenstein ring is a quasi $k$-Gorenstein ring. But the converse does not hold in general. In fact, it is well-known that $k$-Gorenstein rings are left-right symmetric, but the quasi $k$-Gorenstein rings are not left-right symmetrica in general.

If $R$ is an Artin algebra, then each $I_i$ is finitely generated, and  ${\rm fd}_R(I_i)={\rm pd}_R(I_i)$. So $R$ is $k$-Gorenstein if ${\rm pd}_R(I_i)\leq i$ for any $0\leq i\leq k-1$ and $R$ is $k$-Gorenstein for all $k$ if  ${\rm pd}_R(I_i)\leq i$ for any $i\geq 0$. Recall that  an Artin algebra $R$ is said to be of dominant dimension greater than or equal to $k$ if ${\rm pd}_R(I_i)=0$ for $0\leq i\leq k-1$. This is a special case of $k$-Gorenstein rings. There  is a famous Nakayama Conjecture, which states that if an Artin algebra $R$ is of dominant dimension $k$ for all $k$, then $R$ is self-injective \cite{N}. Hence the studying of Artin algebra which are $k$-Gorenstein for some $k$ or all $k$ connects with Nakayama Conjecture.

The following lemma seems to be well-known.
\begin{lem}{\label{injective and flat dimensions}} Let $R/S$  be a Frobenius extension and $M$ a right $S$-module. Then we have the following

$(1)$  r.${\rm id}_R(M\ot_SR)\leq $r.${\rm id}_S(M) $;

$(2)$ r.${\rm fd}_R(M\ot_SR)\leq $r.${\rm fd}_S(M)$.
\end{lem}
\noindent{\bf Proof.} $(1)$ It is trivial that the case of $r.{\rm id}_S(M)=\infty$.

Now, we assume that $r.{\rm id}_S(M)=n$. By the adjoint isomorphism, for any right $R$-module $C$, we have
\begin{align*}
{\rm Ext}^{n+1}_{R^{\rm op}}(C,M\otimes_SR_R)
&\cong{\rm Ext}^{n+1}_{R^{\rm op}}(C, {\rm Hom}_S({_RR}, M))\\
&\cong {\rm Ext}^{n+1}_{S^{\rm op}}(C\ot_RR, M_S)\\
&\cong {\rm Ext}^{n+1}_{S^{\rm op}}(C_S, M_S).
\end{align*}
Then ${\rm Ext}^{n+1}_{R^{\rm op}}(C,M\otimes_SR_R)\cong {\rm Ext}^{n+1}_{S^{\rm op}}(C_S, M_S)=0$. Therefore, $r.{\rm id}_R(M\ot_SR)\leq n$.

$(2)$ Similar to the proof of $(1)$. If $r.{\rm fd}_S(M)=\infty$, it is trivial.

By the associativity of tensor, for any left $R$-module $C$,
we have
\begin{align*}
{\rm Tor}_{n+1}^R(M\otimes_SR_R, C)
&\cong{\rm Tor}_{n+1}^S(M_S, {_SR\ot_RC})\\
&\cong {\rm Tor}_{n+1}^S(M_S, {_SC}).
\end{align*}
If $r.{\rm fd}_S(M)=n$, then ${\rm Tor}_{n+1}^R(M\otimes_SR_R, C)\cong{\rm Tor}_{n+1}^S(M_S, {_SC})=0$. Hence $r.{\rm fd}_R(M\ot_SR)\leq n$.
\hfill$\square$

\begin{prop}{\label{Prop5.3}} Let $R/S$  be a Frobenius extension and $k$ a positive integer. If $S$ is a $k$-Gorenstein ring, then so is $R$.
\end{prop}
\noindent{\bf Proof.} Let $0\rightarrow {S_S}\rightarrow E'_0\rightarrow E'_1\rightarrow\cdots\rightarrow E'_i\rightarrow\cdots$ be a minimal injective resolution of $S$ as a right $S$-module. By assumption, we have $r.{\rm fd}_S(E'_i)\leq i$ for $0\leq i\leq k-1$.
Since $_SR$ is finitely generated projective, we get the following exact sequence
\begin{equation}
\xymatrix@C=0.5cm{0 \ar[r] & {S\ot_SR_R\cong R_R}\ar[r] &{E'_0\ot_SR}\ar[r] & {E'_1\ot_SR} \ar[r] & {\cdots}\ar[r] &{E'_i\ot_SR}\ar[r] &\cdots} \tag{*}
\end{equation}
in mod-$R$. By Lemma \ref{injective and flat dimensions}, the exact sequence $(*)$ is an injective resolution of $R_R$ as a right $R$-module and $r.{\rm fd}_R(E'_i\ot_SR)\leq i$ for $0\leq i\leq k-1$.

Suppose that $0\rightarrow {R_R}\rightarrow I_0^\prime\rightarrow I_1^\prime\rightarrow\cdots\rightarrow I_i^\prime\rightarrow\cdots$ is a minimal injective resolution of $R$ as a right $R$-module. Then we have $r.{\rm fd}_R(I_i^\prime)\leq r.{\rm fd}_R(E'_i\ot_SR)\leq i$ for $0\leq i\leq k-1$ since $I_i^\prime$ isomorphic to a summand of $E'_i\ot_SR$. Hence,$R$ is a $k$-Gorenstein ring.  \hfill$\square$

\begin{cor} Let $R/S$  be a Frobenius extension. If $S$ is an Auslander ring, then so is $R$.
\end{cor}

The fundamental theorem in \cite{B} states that a commutative Noetherian ring $R$ is Iwanaga-Gorenstein (that is, the injective dimension of $R$ is finite) if and only if the flat dimension of the $i$-th term in a minimal injective coresolution of $R$ as an $R$-module is at most $i-1$ for any $i>1$. In the non-commutative case, Auslander and Reiten \cite{AR1} conjectured that, for an Artin algera $R$, $R$ is an Auslander ring implies that  $R$ is Iwanaga Gorenstein.

In \cite{CR}, Chen and Ren shows that $R$ is $k$-Iwanaga Gorenstein if and only if so is $S$ under a Frobenius extension $R/S$ with $R$ is a generator as a left and right $S$-module, respectively. Base on it, there are some natural questions as follows.

\noindent{\bf Questions:} Is the base ring is (quasi) $k$-Gorenstein equivalent to so is the extension ring  under a Frobenius extension $R/S$ with $R$ is a generator as a left and right $S$-module, respectively?
Furthermore, Is the base ring is an Auslander ring equivalent to the extension ring is an Auslander ring under  a Frobenius extension $R/S$ with $R$ is a generator as a left and right $S$-module, respectively?

\section{modules of G-dimension zero}

In this section, we will show that a module over the extension ring is of G-dimension zero if and only if so is its underlying module over the base ring under Frobenius extensions. We point out that we get this assertion, as another application of the results in Section 3, in a completely different way from that in \cite{CH, R} and \cite{Z}.

As the origin of Gorenstein homological algebra, the module of G-dimension zero was introduced in \cite{AB}.
\begin{Def} A module $M$ in $R$-${\rm mod}$ is said to be of G-dimension zero, denoted by G-dim$_R(M)=0$, if it is satisfies:
$(1)$ ${\rm Ext}^i_R(M, R)=0$ for all $i>0$; $(2)$ ${\rm Ext}^i_{R^{\rm op}}({\rm Tr}M, R)=0$ for all $i>0$. We denote by $\mathcal{G}$ the subcategory of $R$-${\rm mod}$ consisting of all modules of G-dimension zero.
\end{Def}

We denote the subcategory of $R$-${\rm mod}$ consisting of all modules $M$ with ${\rm Ext}^i_R(M, R)=0$ for all $i>0$ by ${^\bot R}$, and it is called \emph{the left orthogonal class} of $R$. By definition of $k$-torsionfreeness and the fact that $M$ and ${\rm Tr(Tr}M)$ are projectively equivalent, we have the following observation.

 \begin{lem}{\label{lemma of G-dimension zero}} For an $R$-module $M$, the following are equivalent.

 $(1)$ \rm {G-dim}$_R(M)=0$;

 $(2)$ $M\in{^\bot R}$ and $M$ is $\infty$-torsionfree;

 $(3)$ $M$ and ${\rm Tr}M$ (as a right $R$-module) are both $\infty$-torsionfree.
\end{lem}
\noindent{\bf Proof.} For $(3)$, it is equivalent to the conditions that ${\rm Ext}^i_{R^{\rm op}}({\rm Tr}M, R)=0$ and ${\rm Ext}^i_{R}({\rm Tr(Tr}M), R)=0$ for all integer $i\geq 1$.
Since $M$ and ${\rm Tr(Tr}M)$ are projectively equivalent, the last conditions are equivalent to the conditions that ${\rm Ext}^i_{R^{\rm op}}({\rm Tr}M, R)=0$ and ${\rm Ext}^i_{R}(M, R)=0$ for all integer $i\geq 1$. Thus, we have $(2)\Leftrightarrow (3)$.

As a generalization of module of G-dimension zero, Enochs and Jenda introduced  the Gorenstein projective module for arbitrary modules over a general ring in \cite{EJ}. It follows from \cite[Theorem 4.2.6]{Chris} that a finitely generated module over two-sided Noetherian ring is Gorenstein projective if and only if it is a module of G-dimension zero. There are several different terminologies in the literature for these modules, such as totally reflexive modules, maximal Cohen-Macaulay modules.

  The following lemma shows that the orthogonal classes can be transfered under Frobenius extensions.
\begin{lem}{\label{orthonogal equivalence I}}Let $R/S$  be a Frobenius extension and $M$ an $R$-module. Then $M\in{^\bot R}$ as an $R$-module if and only if $M\in{^\bot S}$ as an $S$-module.
\end{lem}
\noindent{\bf Proof.} Since $R/S$ is a Frobenius extension, we have $\mathbb{T}({_SM})={_RR}\otimes_SM\cong {\rm Hom}_S({_SR_R,{_SM}})=\mathbb{H}({_SM})$ for any
${_SM}\in {\rm mod}S$.
By the adjoint isomorphism, for any $i\geq 0$
\begin{align*}{\rm Ext}^i_S({_SM},S)
&\cong{\rm Ext}^{i}_S({_SR\otimes_RM},S)\\
&\cong{\rm Ext}^{i}_R({_RM}, {\rm Hom}_S({_SR_R}, S))\\
&\cong {\rm Ext}^i_R({_RM}, {_RR\otimes_SS})\\
&\cong {\rm Ext}^i_R({_RM}, {_RR}).
\end{align*}
Consequently, ${\rm Ext}^i_S({_SM},S)=0$ for any $i\geq 1$ if and only if ${\rm Ext}^i_R({_RM}, R)=0$ for any $i\geq 1$. \hfill$\square$

Chen proves in \cite{CH} that the total reflexivity of modules is preserved under the totally reflexive extension, where the totally reflexive extension is a generalization of Frobenius extensions.  Ren and Zhao obtained in \cite{R} and \cite{Z} that the Gorenstein projectivity of modules is preserved under Frobenius extensions by the method of construction of projective coresolution, respectively.
Combining Corollary \ref{torsionfree equ II} and Lemma \ref{orthonogal equivalence I}, in a completely different way from that in \cite{CH, R} and \cite{Z},  we have the following known result, which is due to \cite[Theorem 3.6]{CH}, \cite[Theorem 2.2]{R} or \cite[Theorem 3.2]{Z}.
\begin{thm} {\label{equivalent of G-dimension}}Let $R/S$  be a Frobenius extension and $M$ an $R$-module. Then $M$ is an $R$-module of G-dimension zero if  and only if $M$ is an $S$-module of G-dimension zero.
\end{thm}
\noindent{\bf Proof.} By Lemma \ref{lemma of G-dimension zero}, $M$ is of G-dimension zero as an $R$-module if and only if $M\in{^\bot R}$ and $M$ is $\infty$-torsionfree, which equivalent to $M\in{^\bot S}$ and $M$ is $\infty$-torsionfree as an $S$-module by Lemma \ref{orthonogal equivalence I} and  Corollary \ref{torsionfree equ II}, respectively. And the last conditions equivalent to that $M$ is of G-dimension zero as an $S$-module. \hfill$\square$

\noindent{\bf Acknowledgements}

The author thanks Professor Xiao-Wu Chen for his helpful suggestions. This work is supported by the National Natural Science Foundation of China (No.12371015).
The author thanks the referee for very useful suggestions and comments.

\vspace{0.5cm}




\begin{thebibliography}{50}
\bibitem{A} M. Auslander, Coherent functors, Proceeding Conference on Categorical Algebra, La Jolla 1965, Springer-Verlag, Berlin-Heidelberg-New York (1966), 189-231.
\bibitem{AB} M. Auslander and M. Bridger, Stable module theory, Memoirs Amer. Math. Soc. 94, Amer. Math. Soc., Providence, RI, 1969.
\bibitem{AR2}M. Auslander and I. Reiten, Homologically finite subcategories, in``Proceedings, Tsukuba Workshop on Representations of Algebras", London Math. Soc. Lecture Notes, Vol. 168, Cambridge Univ. Press, Cambridge, UK.
\bibitem{AR}M. Auslander and I. Reiten, $k$-Gorenstein algebras and syzygy modules, J. Pure Appl. Algebra 92 (1994), 1-27.
\bibitem{AR1}M. Auslander and I. Reiten, Syzygy modules for Noetherian rings, J. Algebra 183 (1996), 167-185.
\bibitem{AS} M. Auslander and S. O. Smal${\o}$, Preprojective modules over artin algebras, J. Algebra 66 (1980), 61-122.
\bibitem{B} H. Bass, On the ubiquity of Gorenstein rings, Math. Z. 82 (1963), 8-28
\bibitem{BGS} K. A. Brown, I. G. Gordon and C. H. Stroppel, Cherednik, Hecke and quantum algebras as free Frobenius and Calabi-Yau extensions, J. Algebra 319 (2008), 1007-1034.
\bibitem{CH} X. -W. Chen, Totally reflexive extensions and modules, J. Algebra 379 (2013), 322-332.
\bibitem{CR} X. -W. Chen and W. Ren, Frobenius functors and Gorenstein homological properties, J. Algebra 610 (2022), 18-37.
\bibitem{Chris}L. W. Christensen, Gorenstein dimenensions, Lecture Notes in Math., Vol 1747, Springer-Verlag, Berlin, 2000.
\bibitem{EJ} E. E. Enochs and O. M. G. Jenda,  Gorenstein injective and projective modules, Math. Z. 220 (1995), 611-633.
\bibitem{FGR} R. M. Fossum,  P. A. Griffith and  I. Reiten, Trivial Extensions of Abelian Categories, in: Homological Algebra of Trivial Extensions of Abelian Categories
with Applications to Ring Theory, in: Lecture Notes in Math., vol. 456, Springer-Verlag, Berlin, New York, 1975.
\bibitem{FMS} D. Fischman, S. Montgomery and H.-J. Schneider, Frobenius extensions of subalgebras of Hopf algebras, Trans. Amer. Math. Soc. 349 (1997), 4857-4895.
\bibitem{GT} J. G$\acute{o}$mez-Torrecillas and B. Torrecillas, FTF Rings and Frobenius Extensions, J. Algebra 248 (2002), 1-14.
\bibitem{HH} C. -H. Huang and Z. -Y. Huang, Torsionfree dimension of modules and self-injective dimension of rings, Osaka J. Math. 49 (2012), 21-35.
\bibitem{Hu1} Z. -Y. Huang, Extension closure of $k$-torsionfree modules, Comm. Algebra 27 (1999), 1457-1464.
\bibitem{Hu2} Z. -Y. Huang, $\omega^t$-approximation representations over quasi $k$-Gorenstein algebras, Sci. China (Series A), 42 (1999), 945-956.
\bibitem{HI} Z. -Y. Huang and O. Iyama, Auslander-type conditions and cotorsion pairs, J. Algebra 318 (2007), 93-110.
\bibitem{Kad0} L. Kadison, The jones polynomial and certain separable Frobenius extensions, J. Algebra 186 (1996), 461-475.
\bibitem{Kad} L. Kadison, New examples of Frobenius extensions, University Lecture Series, Vol 14, AMS. Provedence, Rhode Island, 1999.
\bibitem{Kas1} F. Kasch, Grundlagen einer theorie der Frobenius-Erweiterungen, Math. Ann., 127 (1954), 453-474.
\bibitem{N} T. Nakayama, On algebras with complete homology, Abh. Math. Sem. Univ. Hamburg 22 (1958), 300-307.
\bibitem{NT} T. Nakayama and T. Tsuzuku, On Frobenius extension I, Nagoya Math. J. 17 (1960), 89-110; On Frobenius extensions II, Nagoya Math J. 19 (1961), 127-148.
\bibitem{R} W. Ren, Gorenstein projective and injective dimensions over Frobenius extensions, Comm. Algebra 46 (2008), 1-7.
\bibitem{Sch} H. -J. Schneider, Normal basis and transitivity of crossed products for Hopf algebras, J. Algebra 151 (1992), 289-312.
\bibitem{Z} Z. -B. Zhao, Gorenstein homological invariant properties under Frobenius extensions, Sci. China Math. 62 (2019), 2487-2496.

\end{thebibliography}
\end{document}